\documentclass{article}
\usepackage{amssymb}
\usepackage{amsmath}
\usepackage{latexsym}
\usepackage{mathrsfs}
\usepackage{graphicx}
\usepackage[all]{xy}
\newtheorem{Th}{Theorem}
\newtheorem{Lemma}[Th]{Lemma}
\newtheorem{Coro}[Th]{Corollary}
\def\SO{\operatorname{SO}}
\def\SU{\operatorname{SU}}
\def\SL{\operatorname{SL}}

\title{Direct proof of Mckay Correspondence and the representations of finite subgroups of $\SO(4)$}
\author{XIONG Rui}
\begin{document}

\maketitle

\begin{abstract}
  The classic Mckay correspondence gives a connection between finite subgroups of $\SU(2)$ and the simply-laced Dynkin diagrams.
  In this article, a direct proof is presented. The bipartite structure of the Mckay diagrams is introduced.
  After that, the similar method can be used on finite subgroups of $\SO(4)$, we get a edges-coloured graph.
  We finally get some applications about the dimension restriction of the irreducible representation of finite subgroup of $\SO(4)$.
\end{abstract}

\section{Introduction}

In classic group theory, there is a well-known result that the full list of finite subgroups of $\SO(3)$ is
$$C_n,\quad D_{2n}^*,\quad \mathfrak{A}_4, \quad \mathfrak{S}_4,\quad \mathfrak{A}_5, $$
which are cyclic groups, dihedral groups, tetrahedral group, octahedral group, and dodecahedral group.

There is a $2$ to $1$ surjective group homomorphism $\pi:\SU(2)\to \SO(3)$.
Since $-1\in \SU(2)$ is the only order two element, so it is not difficult to find the full list of finite subgroups of $\SU(2)$ is
$$C_n,\quad D_{2n}^*, \quad \mathfrak{A}_4^*, \quad \mathfrak{S}_4^*, \quad \mathfrak{A}_5^*$$
which are cyclic groups, binary dihedral groups, binary tetrahedral groups, binary octahedral groups and binary dodecahedral groups.

Since all the complex representation is unitary, so the above is also exactly the list of all finite subgroups of $\SL_2(\mathbb{C})$.

\def\oto{\makebox[1pc][c]{$-\!\!\!-$}}
\newcommand{\ovto}[1][\circ]{\stackrel{\mbox{$\stackrel{\mbox{$#1$}}\mid$}}}
Mckay correspondence claims that there is a correspondence of them with simply-laced Dynkin diagrams through Mckay diagram as following.
$$\begin{array}{ccc}
C_n & A_{n} & \circ \oto \circ \oto \cdots \oto\circ \oto \circ  \\ [2ex]
D_{2n}^* & D_{n+2} & \circ \oto \circ \oto \cdots \oto \ovto{\circ} \oto \circ \\ [2ex]
\mathfrak{A}_4^* & E_6 & \circ \oto \circ \oto \ovto\circ \oto \circ \oto \circ \\ [2ex]
\mathfrak{S}_4^* & E_7 & \circ \oto \circ \oto \ovto\circ \oto \circ \oto \circ \oto \circ \\ [2ex]
\mathfrak{A}_5^* & E_8 & \circ \oto \circ \oto \ovto\circ \oto \circ \oto \circ \oto \circ\oto \circ \\ [2ex]
\end{array}$$

The original paper is due to Mckay \cite{mckay1980graphs}, but now is not available.
The main method is via the classification of Du Val singularities for example \cite{stekolshchik2008notes}.

Here I will present a pure elementary proof. Finally, I will discuss the representation of the finite subgroups in $\SO(3)$ and $\SO(4)$.

\section{The direct proof}

Before the proof, we need a description of simply-laced Dynkin diagram.
Let $G=(V,E)$ be a finite graph (undirected and possible loops and self-loops).
Without loss of generality assume $V=\{1,\ldots,r\}$. Let us define its Cartan matrix
$(c_{ij})=(2\delta_{ij}-n_{ij})$ where $n_{ij}$ is the number of edges between $i$ and $j$, and $\delta_{ij}$ Kronecker's delta.
\begin{itemize}
  \item If the graph $G$ is connected, and its Cartan matrix is positive definite, then $G$ is simply-laced Dynkin diagram.
Such matrix is positive definite if and only if there exists a vector $(x_{j})$ with $x_j>0$, such that $\sum_{j=1}^r c_{ij}x_j>0$ for all $i$.

  \item If the graph $G$ is connected, and its Cartan matrix is positive semi-definite, then $G$ is simply-laced Euclidean diagram (or, extended Dynkin diagram).
Such matrix is positive semi-definite if and only if there exists a vector $(x_{j})$ with $x_j>0$, such that $\sum_{j=1}^r c_{ij}x_j=0$ for all $i$.

  \item The full list of simply-laced Euclidean diagrams is
$$\begin{array}{cc}
\tilde{A}_n & \xymatrix@R=0pc{& &&& \\&\circ \ar@{-}[r] & \cdots \ar@{-}[r] & \circ \ar@{-}`[ur]`[u]`[lll]|{\mbox{$\circ$}}`[ll][ll]&}\\[2ex]
\tilde{D}_n & \circ \oto \ovto\circ \oto \cdots \oto  \ovto\circ \oto \circ \\[2ex]
\tilde{E}_6 &  \circ \oto \circ \oto \ovto{\ovto\circ} \oto \circ \oto \circ \\ [2ex]
\tilde{E}_7 &  \circ \oto\circ \oto \circ \oto \ovto\circ \oto \circ \oto \circ \oto \circ \\ [2ex]
\tilde{E}_8 &  \circ \oto \circ \oto \ovto\circ \oto \circ \oto \circ \oto \circ \oto \circ \oto \circ \\ [2ex]
\end{array}$$
\end{itemize}

Above facts can be found for example in \cite{benson1991representations} volume 1 page 120 theorem 4.5.8.
Note that in the book, a more general diagram (not necessary undirected) is defined.

\bigbreak
Now, let us turn to representation theory.
\def\bbC{{\mathbb{C}}}
\def\Hom{\operatorname{Hom}}
Generally, let us consider an $n$-dimensional representation $W$ of finite group $G$ over $\bbC$.
Let $\{V_i\}_{i=1}^r$ be the full list of pairwise
non-isomorphic irreducible representations of $G$. Let
$$n_{ij}=\textrm{The multiplicity of $V_i$ in $W\otimes_{\bbC} V_j$}=\Hom_G(V_i,W\otimes_k V_j). $$
We define a graph with vertices $\{1,\ldots, r\}$, and connected $i$ and $j$ by $n_{ij}$ arrows.
This is known as \emph{Mckay diagram} with respect to $W$.

Here are some properties.
\begin{itemize}
  \item If $W$ is self-dual, then the graph is undirected. That is, if $W^\vee\cong W$ as $G$-representations, then $n_{ij}=n_{ji}$.
  Since
  $$\dim\Hom_G(V_i,W\otimes_k V_j)=\dim\Hom_G(V_i\otimes_k W^{\vee}, V_j)=\dim \Hom_G(V_j,W\otimes_k V_i). $$
  \item If $W$ is faithful, then the graph is connected. Since we know that any irreducible representation appear in some fold of tensor product of faithful representation, that is, $\Hom_G(V_i,W^{\otimes n})\neq 0$ for some $n$.
      More precisely, if it is not connected, then we can decompose $\{V_i\}_{i=1}^r$ by $\mathfrak{V}_1\sqcup \mathfrak{V}_2$.
      Then $W^{\otimes n} \otimes V_i$ decomposes into irreducible representations in $\mathfrak{V}_\bullet$
       for any $V_i\in \mathfrak{V}_\bullet$. Now, consider $V_i\otimes V_j$ for $V_i\in \mathfrak{V}_1$ and $V_j\in \mathfrak{V}_2$ a contradiction.
  \item The vector $(\dim V_j)_j$ is an eigenvector belonging to $\dim W$ of $(n_{ij})_{i,j}$, since
  $$\begin{array}{rl}
  \sum_{j=1}^r n_{ij}\dim V_j & = \sum_{j=1}^n \dim \Hom_G(V_i,W\otimes V_j)\dim V_j \\
  & =  \dim \Hom_G\big(V_i,W\otimes \bigoplus_{j=1}^n (\dim V_j)V_j\big)\\
  & = \dim  \Hom_G\big(V_i,W\otimes \bbC[G]\big)\\
  & = \dim \Hom_{\bbC}(V_i,W)\\
  & = \dim W\cdot \dim V_i.
  \end{array}$$
\end{itemize}

\bigbreak

Our situation is when $n=2$, since now $G\subseteq \SL_2(\mathbb{C})$,
the natural two-dimensional representation $W$ is automatically faithful and self-dual by considering the character.
So by above, the Cartan matrix of the Mckay diagram annulates a positive vector $(\dim V_i)$, so the diagram is Euclidean diagram from the description above.
If we discard the trivial representation, say the \emph{reduced Mckay diagram}, it will be a Dynkin diagram.

We have the following list, with $\times$ the trivial representation, labelled numbers the dimensions of representations.
$$\begin{array}{ccc}
C_n & \tilde{A}_{n-1} & \xymatrix@R=0pc{& &&& \\&1 \ar@{-}[r] & \cdots \ar@{-}[r] & 1 \ar@{-}`[ur]`[u]`[lll]|{\mbox{$\times$}}`[ll][ll]&}  \\ [2ex]
D_{2n}^* & \tilde{D}_{n+2} & 1 \oto \ovto[\times]2 \oto \cdots \oto  \ovto[1]2 \oto 1 \\ [2ex]
\mathfrak{A}_4^* & \tilde{E}_6 & 1 \oto 2 \oto \ovto[\times]{\ovto[2]3} \oto 2 \oto 1  \\ [2ex]
\mathfrak{S}_4^* & \tilde{E}_7 & \times \oto 2 \oto 3 \oto \ovto[2]4 \oto 3 \oto 2 \oto 1 \\ [2ex]
\mathfrak{A}_5^* & \tilde{E}_8 & 2 \oto 4 \oto \ovto[3]6 \oto 5 \oto 4 \oto 3\oto 2\oto \times \\ [2ex]
\end{array}$$

This proves the Mckay correspondence.

\begin{Th}[Mckay correspondenc \cite{mckay1980graphs}]\label{ClassicMckay}If $G$ is a finite subgroup of $\SL_2(\mathbb{C})$, then the reduced Mckay diagram is a simply-laced Dynkin diagram.
\end{Th}

\section{Finite subgroups of $\SO(3)$}

Let $G$ be a finite subgroup of $\SU(2)$, and denote $\pi(G)\subseteq \SO(3)$ the image under the morphism $\pi:\SU(2)\to \SO(3)$.
We assume the order of $G$ is divided by $2$ --- excluding the cases when $G$ is cyclic of odd order.
Now $-1\in G$, and $2|\pi(G)|=|G|$. Now, $\pi(G)$ is the group which we are familiar with.
We will still denote $W$ the natural two-dimensional representation.

Note that
\begin{itemize}
  \item Each (irreducible) representation of $\pi(G)$ is also an (irreducible) representation of $G$ by the natural map $G\to \pi(G)$.
  \item A representation $V$ of $G$ induces a representation of $\pi(G)$ if and only if $-1$ acts trivially on $V$.
\end{itemize}
Denote $\mathfrak{V}\subseteq \{V_i\}_{i=1}^r$ inducing representation of $\pi(G)$,
and $\mathfrak{V}'$ the rest of them.

Let the character of $V$ be $\chi$, since the character is sum of $\dim V$ many roots of unity, so $-1$ acts trivially
if and only if $\chi(-1)=\dim V=\chi(1)$.
Note that the character $\varphi$ of $W$ satisfy $\varphi(-1)=-2=-\varphi(1)$. As a result, we have the following.
\begin{itemize}
  \item If $V_i\in \mathfrak{V}$, whose character $\chi_i$ satisfies $\chi_i(-1)=\pm\chi_i(1)$,
  then each of its neighborhood $V_j$ in Mckay diagram satisfies $\chi_i(-1)=\mp\chi_i(1)$.

  Assume $W\otimes V_i=\bigoplus_{j\in N(i)} n_{ij} V_j$, where $N(i)=\{j: n_{ij}\neq 0\}$ the neighborhood of $i$.
  The decomposition gives rise to $\varphi\chi_i = \sum_{j\in N(i)} n_{ij} \chi_j$.   So
  $$\mp \varphi(1)\chi_i(1)=\varphi(-1)\chi_i(-1)=\sum_{j\in N(i)} n_{ij} \chi_j(-1).$$
  But $\chi_j(-1)$ is sum of $\chi(1)$ many roots of unity, so $|\chi_j(-1)|\leq |\chi_j(1)|$.
  As a result, to achieve the equality, $\chi_j(-1)=\mp \chi_j(1)$.

  \item Since the Mckay diagram is connected, and trivial representation is of course the above case, so all irreducible representation $V_i$ having its character $\chi_i$ satisfy $\chi_i(-1)=\pm \chi_i(1)$.
  \item The neighborhoods of representation form $\mathfrak{V}$ are all from $\mathfrak{V}'$.
  Et vice versa, the neighborhood of representation from $\mathfrak{V}'$ are all from $\mathfrak{V}$.
\end{itemize}
The result is summarized in the following theorem.

\begin{Th}\label{BipartitegraphEuclidean}
The representation from $\pi(G)$ and the representation not from divide the Mckay diagram into a bipartite graph.
\end{Th}

We can read the dimensions from the diagram.
$$\xymatrix@R=0pc{
1 \ar@{-}[r] & 1 \ar@{-}[dl]\\
1 \ar@{-}[r] & \vdots\ar@{--}[dl] \\
\vdots \ar@{-}[r]& 1 \ar@{-}[dl]\\
1 \ar@{-}[r] & 1\ar@{-}[luuu]}\quad
\xymatrix@R=0pc{
1\ar@{-}[dr] \\
1 \ar@{-}[r] & 2 \ar@{-}[dl]\\
2 \ar@{-}[r] & \vdots \ar@{--}[dl]\\
\vdots \ar@{-}[r]& 2 \ar@{-}[dl]\\
2 \ar@{-}[r]\ar@{-}[dr] & 1\\
& 1}\quad
\xymatrix@R=0pc{
1\ar@{-}[dr] \\
1 \ar@{-}[r] & 2 \ar@{-}[dl]\\
2 \ar@{-}[r] & \vdots \ar@{--}[dl]\\
\vdots \ar@{-}[r]& 2 \ar@{-}[dl]\ar@{-}[ddl]\\
1 & \\
1 & }\quad
\xymatrix@R=0pc{
1\ar@{-}[r] & 2 \\
1\ar@{-}[r] & 2 \\
1\ar@{-}[r] & 2\\
3\ar@{-}[uuur]\ar@{-}[uur]\ar@{-}[ur]}\quad
\xymatrix@R=0pc{
1\ar@{-}[dr] &  \\
3\ar@{-}[r]\ar@{-}[dr] & 2 \\
2\ar@{-}[r] & 4\\
3\ar@{-}[r]\ar@{-}[ur]& 2\\
1\ar@{-}[ur]}\quad
\xymatrix@R=0pc{
1\ar@{-}[dr]\\
3\ar@{-}[dr]\ar@{-}[r] & 2 \\
5\ar@{-}[r]\ar@{-}[dr] & 4 \\
3\ar@{-}[r] & 6\ar@{-}[dl]\\
4\ar@{-}[r]& 2}\eqno{(*)}$$
Or more beautiful picture,
$$\begin{array}{cccc}
G & \pi(G) & & \\
C_{2k}& C_k & \tilde{A}_{2k} & \xymatrix@R=0pc{& &&&&& \\&\circ \ar@{-}[r] &\bullet \ar@{-}[r] & \circ \ar@{-}[r] & \cdots \ar@{-}[r] & \circ \ar@{-}`[ur]`[u]`[lllll]|{\mbox{$\bullet$}}`[llll][llll]&}  \\ [2ex]
D_{4k}^* & D_{2k}& \tilde{D}_{2k+2}  & \bullet \oto \ovto[\bullet]\circ \oto \bullet \oto \circ \oto \cdots \oto  \ovto[\circ]\bullet \oto \circ \\ [2ex]
D_{4k+2}^* & D_{4k+2} &\tilde{D}_{2k+3}  & \bullet \oto \ovto[\bullet]\circ \oto \bullet \oto \circ \oto \cdots \oto  \ovto[\bullet]\circ \oto \bullet \\ [2ex]
\mathfrak{A}_4^*&\mathfrak{A}_4 & \tilde{E}_6 & \bullet \oto \circ \oto \ovto[\bullet]{\ovto[\circ]\bullet} \oto \circ \oto \bullet  \\ [2ex]
\mathfrak{S}_4^*& \mathfrak{S}_4 & \tilde{E}_7 & \bullet \oto \circ \oto \bullet \oto \ovto[\bullet]\circ \oto \bullet \oto \circ \oto \bullet \\ [2ex]
\mathfrak{A}_5^*& \mathfrak{A}_5 & \tilde{E}_8 & \circ \oto \bullet \oto \ovto[\bullet]\circ \oto \bullet \oto \circ \oto\bullet\oto \circ\oto \bullet\\ [2ex]
\end{array}$$
This reflects the dimensions of the the representation of finite subgroups of $\SO(3)$.
The presentation of $D_{2n}$ can be found in \cite{serre1977linear} page 37 (the notation is $D_n$ for $D_{2n}$ here).
The presentation of $\mathfrak{A}_4,\mathfrak{S}_4,\mathfrak{A}_5$ can be found in \cite{fulton2013representation}
page 18, page 19 and page 29 respectively.

\section{Finite subgroups of $\SO(4)$}

We know the universal cover of $\SO(4)$ is exactly $\SU(2)\times \SU(2)$.
Let $G$ be a finite subgroup of $\SU(2)\times \SU(2)\subseteq \SL_4(\mathbb{C})$, and $\pi(G)\subseteq \SO(4)$ the image under covering map.
Assume as well, $-1\in G$, that is, $|G|=2|\pi(G)|$.
Let $W$ be the natural representation of dimension $4$.
Let $G_1$ the image under the projection of the first $\SL_2(\mathbb{C})$, and $G_2$ the second.
Let $W_i$ be the two-dimensional representation of $G_i$ for $i=1,2$.
It is easy to see $W\cong W_1\oplus W_2$.
Let $\{V_i\}_{i=0}^r$ be the full list of pairwise nonisomorphic irreducible representations.
Denote
$$n_{ij}^k=\dim\Hom(V_i,W_k\otimes V_j), \quad k=1,2,\varnothing. $$
So $n_{ij}=n_{ij}^1+n_{ij}^2$.

Denote the Mckay diagram of $W$ to be $\Gamma$.
We can colour the edges between $i$ and $j$ by $n_{ij}^1$ many $1$'s,
and $n_{ij}^2$ many $2$'s. Let us denote $\Gamma_i$ be the subgraph of all vertices and all edges coloured by $i$ for $i=1,2$.
Similar to what we did last sections, we have the following properties.
\begin{itemize}
  \item Since $W_k$ is also faithful and self-dual, $\Gamma_k$ is undirected and connected for $k=1,2,\varnothing$, that is, $n_{ij}^k=n_{ji}^k$.
  \item By considering the connected component of $\Gamma_i$, $\Gamma_i$ is disjoint union of Euclidean diagram for $i=1,2$.
  \item As what we did for $\SO(3)$ in theorem \ref{BipartitegraphEuclidean}, the diagram is also bipartite.
\end{itemize}

\begin{Lemma}\label{intersectingtraversally}
The Euclidean diagram from $\Gamma_1$ and $\Gamma_2$ intersect transversally.
That is, each connected component $E_1\subseteq \Gamma_1$ and $E_2\subseteq \Gamma_2$ intersect.
\end{Lemma}

Pick $V_i\in E_i$. Note that $E_i$ is exactly the irreducible representations appearing in $W_i^{\otimes n}\otimes E_i$ for some $n$.
So we want to show that $\Hom_G(W_i^{\otimes n}\otimes E_i, W_j^{\otimes m}\otimes E_j)\neq 0$ for some $m,n$.
It suffices to show $W_i^n\otimes W_j^{m}$ is faithful for some $m,n$. If some $g\in G$ such that
$\big(\frac{\varphi_1(g)}{\varphi_1(1)}\big)^n\big(\frac{\varphi_2(g)}{\varphi_2(1)}\big)^m=1$ for some $m,n$,
where $\varphi_i$ is the character of $W_i$ for $i=1,2$.
Then $\big|\frac{\varphi_1(g)}{\varphi_1(1)}\big|=\big|\frac{\varphi_2(g)}{\varphi_2(1)}\big|=1$
which implies $g$ acts as scalar over $W_i$ for $i=1,2$. Since $G$ acts over $W_i$ through $\SL_2(\mathbb{C})$,
it is only possible when $g=\pm 1$. So taking $W_1\otimes W_2 \otimes W_2$ works.

\begin{Th}Assume the Mckay diagram $\Gamma=(V,E)$,
then $\Gamma$ is bipartite and $E$ admits a decomposition $E=E_1\sqcup E_2$, such that
\begin{itemize}
  \item $\Gamma_i=(V,E_i)$ is disjoint union of Euclidean diagrams for $i=1,2$.
  \item The Euclidean diagrams from $\Gamma_1$ and $\Gamma_2$ intersect transversally.
\end{itemize}
\end{Th}

\begin{Th}The vector $(\dim V_j)$ is the only vector $(x_j)$ up to scalar satisfying $\sum n_{ij}x_j =4x_i$ for any $i$.
\end{Th}

Firstly, the spectral radium of the adjacency matrix of Euclidean diagram is $2$.
By Frobenius-Perron's argument, see \cite{etingof2016tensor} page 51 Theorem 3.2.1, there is a positive eigenvector  $(x_{j})$ with $x_j>0$ belonging to the maximal eigenvalue $\lambda$.
If $\lambda>2$, then $\sum_{j} 2\delta_{ij}-n_{ij}x_j=(2-\lambda)x_i>0$, which implies $(2\delta_{ij}-n_{ij})$ is positive definite, it is impossible.

Since the adjacency matrix is symmetric, so the eigenvalues are exactly the singular values so
$$\sum_{i}\big|\sum_j n_{ij} x_j\big|^2\leq \sum_{j}|2x_j|^2.$$

Assume $(x_j)$ satisfying $\sum n_{ij}x_j =4x_i$ for any $i$. Considering the modulus,
we find over each Euclidean diagram $E$ (coloured by $1$ or $2$), $(x_j)_{j\in E}$ is determined up to a scalar.
But the diagram is connected, so $(x_j)$ is determined up a scalar.

\begin{Coro}\label{EachEuclideanDiagram}In the proof above, we proved that over each Euclidean diagram $E\subseteq G_{1,2}$
such $(x_j)$ is a scalar of $(*)$ after theorem \ref{BipartitegraphEuclidean}.
\end{Coro}

\begin{Coro}\label{Readtheorder}The order of group is determined by the Mckay diagram.
\end{Coro}

Since we can pick a $(x_j)\neq 0$ such that $\sum n_{ij}x_j =4x_i$ for any $i$.
Now $x_i>0$, we can assume the minimal $x_i$ is $1$, now $x_i=\dim V_i$, so the order of group is $\sum_{i=1}^r x_i^2$.

%
%
%
%

\section{Examples}

A direct example is the product of two finite subgroups of $\SU(2)$,
assume $G=G_1\times G_2$. Let $\{U_{ik}\}_i$ is the full list of irreducible representations of $G_k$ for $k=1,2$, then
$$\begin{cases}
\dim\Hom_G(V_{i1}\otimes V_{\imath2}, W_1\otimes V_{j1}\otimes V_{\jmath2})=\dim \Hom_G(V_{i1},W_1\otimes V_{j1})\delta_{\imath\jmath}\\
\dim\Hom_G(V_{i1}\otimes V_{\imath2}, W_2\otimes V_{j1}\otimes V_{\jmath2})=\dim \Hom_G(V_{\imath 1},W_2\otimes V_{\jmath1})\delta_{ij}
\end{cases}$$
then the edge-coloured Mckay diagram is a product of two edge-coloured Euclidean diagrams.
For example the following.
$$\includegraphics[width=0.25\linewidth]{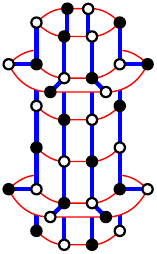}$$
Actaully, we have the following.

\begin{Th}The finite subgroup of $\SL(2)\times \SL(2)$ is a product of two finite subgroups in $\SL(2)$ (up to an isomorphism) if and only if the edge-coloured Mckay diagram is a product of two edge-coloured Euclidean diagrams.
\end{Th}

Since we can read $G_1$ and $G_2$ from the colour, now $G\subseteq G_1\times G_2$,
and by reading the order of the group by corollary \ref{Readtheorder}, it takes the equality.

\bigbreak
Another example is the diagonal the finite group in the diagonal of $\SU(2)\times \SU(2)$. Assume $G_1=G_2$, and
$G\subseteq G_1\times G_2$ is the diagonal. Now $W_1\cong W_2$, so the Mckay diagram is an Eulidean diagram with edges double.
For example the following.
$$\includegraphics[width=0.5\linewidth]{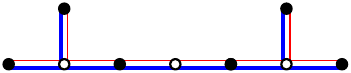}$$
Actaully, we have the following.

\begin{Th}The finite subgroup of $\SL(2)\times \SL(2)$ is a finite subgroup in the diagonal (up to an isomorphism) if and only if the edge-coloured Mckay diagram is an Eulidean diagram with edges double.
\end{Th}

Since we can firstly read $G_1\cong G_2$ from the coloured diagram, by read the dimensions by corollary \ref{Readtheorder}.
Now $G\subseteq G_1\times G_2$ with the projection $G\to G_1$ and $G\to G_2$ isomorphisms, so it is isomorphism to a finite group in the diagonal.



\section{Applications}

Now turn to some applications.
The main tool is corollary \ref{EachEuclideanDiagram} and $(*)$ after theorem \ref{BipartitegraphEuclidean}.
Of course, since the finite subgroups of $\SO(4)$ is classified, see for example \cite{wiki:PointGroupsinFourDim},
so all of the application has a ``violent proof''.

\begin{Th}The irreducible representation of any finite subgroup $G$ of $\SO(4)$ is of dimension no more than $36$.
\end{Th}

Since every representation is connected with trivial representation by at most two Euclidean diagram,
so from $(*)$ after theorem \ref{BipartitegraphEuclidean} and corollary \ref{EachEuclideanDiagram}.
The dimensions are bounded by $36$.


Of course, such bound is achieved, by the image of $\mathfrak{A}_5^*\times \mathfrak{A}_5^* \subseteq \SU(2)\times \SU(2)$.

\begin{Th}If $G$ is a finite subgroup of $\SO(4)$,
then any prime divisor of the dimension of an irreducible representations of $G$ is $2$, $3$ or $5$.
\end{Th}

By corollary \ref{EachEuclideanDiagram}, and over each Euclidean diagram, if $n>2$ is a dimension, then $n/3$, $n/5$ or $n/2$ is a dimension.

\begin{Th}If $G$ is a finite subgroup of $\SO(4)$, then $G$ has an irreducible representation of odd dimension more than one iff and only if it has one of dimension $3$.
\end{Th}

By corollary \ref{EachEuclideanDiagram}, the odd dimension $n$ must divided by $3$ and $5$, and $n/3$ or $n/5$ is also a dimension.
So we can assume that $n=5$, but now, it contains $\tilde{E}_8$ with scalar $1$, so it also has $3$ in dimensions.

Note that $\mathfrak{S}_5$ can be embeded into $\operatorname{O}(4)$ as symmetric group of $5$-cell.
But we have the following result.

\begin{Coro}The symmetric group $\mathfrak{S}_5$ cannot be embedded into $\SO(4)$.
\end{Coro}

Because the dimension of all the irreducible representations is $(1,1,4,4,5,5,6)$ containing $5$ but excluding $3$.


\bibliographystyle{plain}
\bibliography{bibfile}
\vfill

\end{document}